\documentclass{article}

\usepackage{fullpage,amsthm,amsmath,amsfonts,amssymb}

\usepackage{color}

\usepackage{latexsym}
\usepackage{enumerate,dsfont,pifont}
\usepackage{graphicx}

\newtheorem*{Thm}{Theorem}
\newtheorem*{Lem}{Lemma}

\theoremstyle{definition}
\newtheorem*{Rem}{Remark}
\newtheorem*{ack}{Acknowledgement}

\newcommand{\C}{\ensuremath{\mathbb{C}}}

\newcommand{\R}{\ensuremath{\mathbb{R}}}
\newcommand{\rn}{\ensuremath{{\mathbb{R}^\dimension}}}

\newcommand{\N}{\ensuremath{\mathbb{N}}}

\newcommand{\Dskript}{\ensuremath{\mathcal{D}}}

\renewcommand{\leq}{\leqslant}
\renewcommand{\geq}{\geqslant}
\renewcommand{\Re}{\ensuremath{{\rm Re\,}}}

\renewcommand{\phi}{\varphi}

\newcommand{\dimension}{d}

\newcommand{\keywords}[1]{#1}
\newcommand{\ccode}[1]{#1}

\newcommand{\bbc}{\mathbb{C}}
\newcommand{\bbr}{\mathbb{R}}
\newcommand{\bbp}{\mathbb{P}}
\newcommand{\bbe}{\mathbb{E}}

\newcommand{\bbn}{\mathbb{N}}

\newcommand{\cb}{\mathcal{B}}
\newcommand{\cf}{\mathcal{F}}

\newcommand{\cp}{\mathcal{P}}

\newcommand{\abs}[1]{\left| #1 \right|}
\newcommand{\norm}[1]{\left\| #1 \right\|}

\newcommand{\loi}[2]{\left] #1 , #2 \right]}                              
\newcommand{\roi}[2]{\left[ #1 , #2 \right[}

\newcommand{\mutilde}{\widetilde{\mu}}
\DeclareMathOperator{\id}{id}

\newcommand{\pushforward}{{\scriptstyle\bullet}}

\newcommand{\stepsize}{h}

\begin{document}

\title{The Euler scheme for Feller processes}

\author{Bj\"orn B\"ottcher\\ \textit{\small TU Dresden, 
Institut f{\"u}r mathematische Stochastik,}\\
\textit{\small 01062 Dresden,
Germany,
{bjoern.boettcher at tu-dresden.de}}\\[0.5cm]
 Alexander Schnurr\\
\textit{\small TU Dortmund, Fakult\"at f\"ur Mathematik, Vogelpothsweg 87,}\\
\textit{\small 44227 Dortmund, Germany, Alexander.Schnurr at math.tu-dortmund.de}\\
 }

\maketitle

\begin{abstract}

We consider the Euler scheme for stochastic differential equations with jumps, whose intensity might be infinite and the jump structure may depend on the position. This general type of SDE is explicitly given for Feller processes and a general convergence condition is presented.

In particular the characteristic functions of the increments of the Euler scheme are calculated in terms of the symbol of the Feller process in a closed form. These increments are increments of L\'evy processes and thus the Euler scheme can be used for simulation by applying standard techniques from L\'evy processes.

\end{abstract}

\medskip

\ccode{AMS Subject Classification: \emph{Primary}:
60H35, 
\emph{Secondary}:
65C30, 
60J75, 
60J25, 
47G30. 
}

\section{Introduction and main result}

The most general stochastic differential equation (SDE) defining a time homogeneous Markov process $(X_t)_{t\geq 0}$ taking values in $\R^d, d\geq 1$ is of the form
\begin{equation} \label{sde}
\begin{split}
X_t = X_{0} &+ \int_{0}^t a(X_{s-})\ ds + \sum_{j=1}^n \int_{0}^t b^{(j)}(X_{s-})\ dW^{(j)}_s\\
&+ \int_{0}^t \int_{|u|\leq 1} k(X_{s-},u)\ q(\cdot ;ds,du)\\
&+ \int_{0}^t \int_{|u|> 1} k(X_{s-},u)\ p(\cdot ;ds,du)
\end{split}
\end{equation}
where $a,b, k$ are the \textit{coefficients}, $W^{(j)}$ are independent Brownian motions, $p$ is a Poisson random measure and $q$ is the corresponding compensated Poisson random measure  (cf.\,Skorokhod \cite{Skor65}). This equation even includes time inhomogeneous Markov processes, since one can transform any time inhomogeneous Markov process by extending the state space into a time homogeneous Markov process (cf.\,Wentzell \cite{Went79} 8.5.5).

Note that letting $k\equiv 0$ in \eqref{sde} yields a diffusion
equation, the classical setting for the Euler-Maruyama scheme. A
L\'evy driven SDE is also just a special case of this equation. To
see this let $f$ be a $d\times n$ valued function, $l \in \R^n$,
$\sigma$ a positive semi definite matrix in $\R^{n\times n}$ and $N$
an $n$-dimensional L\'evy measure and set $a(X_{s-}) = f(X_{s-}) l,
b(X_{s-}) = f(X_{s-}) \sigma, k(X_{s-},u) = f(X_{s-})$ and let $p$
have the dual predictable projection $ds\,N(du)$. Then equation
$X_t = X_{0} + \int_{0}^t f(X_{s-})dZ_s$ where $Z_s$ is the
L\'evy process on $\R^n$ with triplet $(l,\sigma^2,N).$ In contrast
to these two examples $k$ may depend on $u$ in this note. For an
overview of numerical approximation schemes for this case with
finite jump intensity see for example Bruti-Liberatia and Platen
\cite{BrutPlat2007}.

The Euler approximation with step size $\stepsize$ for an SDE of form \eqref{sde} is given by
$$\bar{X}_0 := X_0 $$
and for $m\in \N_0$
\begin{equation} \label{euler}
\begin{split}
\bar{X}_{(m+1)\cdot \stepsize} := \bar{X}_{m\cdot \stepsize} &+ \int_{m\cdot \stepsize}^{(m+1)\cdot \stepsize}a(\bar{X}_{m\cdot \stepsize})\ ds + \sum_{i=1}^d \int_{m\cdot \stepsize}^{(m+1)\cdot \stepsize} b^{(j)}(\bar{X}_{m\cdot \stepsize})\ dW^{(j)}_s\\
&+ \int_{m\cdot \stepsize}^{(m+1)\cdot \stepsize} \int_{|u|\leq 1} k(\bar{X}_{m\cdot \stepsize},u)\ q(\cdot ;ds,du)\\
&+ \int_{m\cdot \stepsize}^{(m+1)\cdot \stepsize} \int_{|u|> 1} k(\bar{X}_{m\cdot \stepsize},u)\ p(\cdot ;ds,du).
\end{split}
\end{equation}
For the convergence of the Euler scheme it is necessary that small changes of $X_0$ only cause comparable small changes of the distribution of $X_t$ for fixed $t>0$ and for $t \downarrow 0$ the distribution of $X_t$ should converge to the Dirac distribution with point mass at $X_0$. A natural choice of processes satisfying these conditions are Feller processes (cf.\,Section \ref{def}).

Stroock \cite{Stro2003} uses an Euler scheme approach to construct
Feller processes as solutions to \eqref{sde}, although the SDE is
not mentioned explicitly. Therein conditions are formulated in terms
of the coefficients, they are related to the usual Lipschitz
conditions which ensure the existence of a solution to \eqref{sde}
(see \cite{Skor65}).

Contrary to this, the conditions in
the theorem below will be given in terms of the symbol of the generator rather than in
terms of the coefficients of the SDE. This is motivated by the
following facts: 

(i) For construction and analysis of Feller processes the
generator is the natural object to start with, see for example Ethier and Kurtz \cite{EthiKurt86} Chapter 1. Furthermore using formula (3.13) of Courr\`{e}ge \cite{Cour66} it
is possible to
 calculate the symbol without knowledge of the SDE.
Nevertheless, it is a natural question, if the process can be
described by an SDE and then approximated (and simulated) by an
Euler scheme. Our main theorem gives an affirmative answer to this
question. 

(ii) If the coefficients of an
SDE of type \eqref{sde-feller} below are given, one can directly write
down the symbol by formula \eqref{symbolofeq} and check if the
assumptions of the Theorem are fulfilled.

(iii) If the process under consideration is given by a
different type of SDE (cf.\,M\'etivier \cite{Meti1982} Chapter 8) it
is sometimes hard to transform it into the Skorokhod-type. The
symbol on the other hand can occasionally be written down directly
and in a neat way: in \cite{SchiSchn09} it is shown that the
symbol of the solution of the L\'evy driven SDE $X_t = X_{0} + \int_{0}^t f(X_{s-})dZ_s$  is
$\psi(f(x)'\xi)$ where $\psi$ is the characteristic exponent of the
driving L\'evy process $Z_s$.

(iv) While the coefficients depend on the choice of the
SDE-type and the truncation function (in \eqref{euler} we have
chosen $1_{\{\abs{u}\leq 1\}}$; someone interested in limit theorems
would probably choose a continuous function), this is not the case
for the symbol. In this sense the symbol is a `canonical object'.

Thus the conditions in our main theorem are stated in terms of the generator and its symbol: 

\begin{Thm}
Let $(X_t)_{t\geq 0}$ be a Feller process with generator $A$. Assume that
\begin{gather}
    \tag{A1}\label{assumptioncore}
    C_c^\infty(\R^d) \text{\ \ is an operator core of $A$, i.e.\ the closure of }
    {A\big|_{C_c^\infty(\bbr^d)}} \text{ is } A.
\intertext{Let $q(x,\xi)$ be the symbol of $A\big|_{C_c^\infty(\bbr^d)}$ and assume that}
    \tag{A2}\label{assumptionboundedness}
    \exists\, c >0 : |q(x,\xi)|\leq c (1+|\xi|^2) \text{\ \  for all $x$ and $\xi$},\\
    \tag{A3} \label{assumption-nokill}
  q(x,0)=0 \text{\ \ for all $x$}.
\end{gather}
Then the Euler scheme \eqref{euler} for the corresponding SDE converges to $(X_t)_{t\geq 0}$ weakly in $D([0,\infty),\R^d)$, moreover given that $\bar{X}_{(m)\cdot \stepsize} =x$ the next step of the Euler scheme $\bar{X}_{(m+1)\cdot \stepsize}$ has the characteristic function
\begin{equation}\label{char-fun-increment}
e^{ix'\xi} e^{-\stepsize q(x,\xi)}.
\end{equation}
\end{Thm}
Since a Feller process is a time homogeneous Markov process it is the solution of an SDE of the form \eqref{sde}, this SDE is meant by \textit{the corresponding SDE}. The SDE will be explicitly given in Section \ref{proof} and the definition of the other terms and objects appearing in the Theorem can be found in the next section.

\begin{Rem}
For simulations formula \eqref{char-fun-increment} is the key.
It shows that starting at $x$ the next position of the scheme is the
sum of $x$ and the increment (over time $\stepsize$) of a L\'evy
process with characteristic exponent $\xi \mapsto q(x,\xi)$. Thus
the simulation of Feller processes using the Euler scheme is
obvious, if one knows how to simulate L\'evy increments. For the
latter several techniques are well known, see for example Cont and Tankov \cite{ContTank2004}.

An example of a simulation is given in Figure \ref{fig:stable}. It shows a simulated sample path of a one dimensional stable-like process with generator $-(-\Delta)^{\alpha(x)/2}$ where $\alpha(x) = ((0.9+x) \land 1.9) \lor 0.9$, i.e. the process behaves almost like Brownian motion (with double speed) if $X_t > 1$ and almost like a Cauchy process if $X_t< 0$. The state space dependent behavior can be nicely observed in the figure. For the existence of the process and its properties see for example \cite{Bass1988a}.

Further properties of the scheme as speed of convergence and error estimates are part of ongoing research. A practitioners guide to simulation of Feller processes will be given in \cite{B2010}.
\begin{figure}[!ht]
\begin{center}
\includegraphics[width=4in]{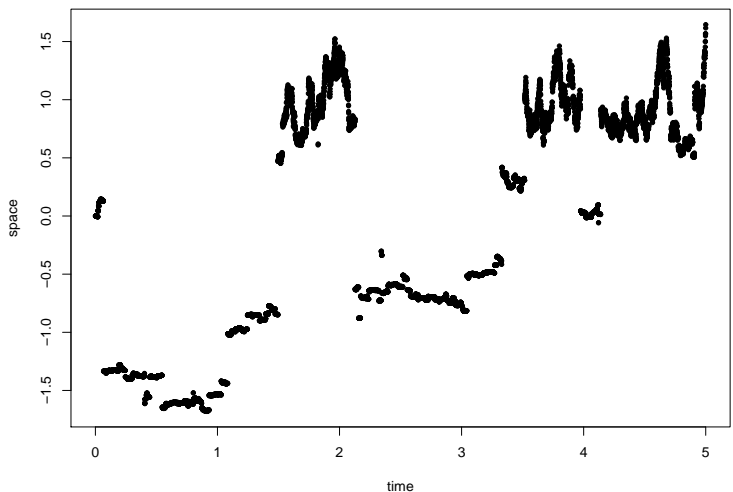} 
\end{center}
	\caption{Stable-like process, $T = 5$ with 1000 steps.}
	\label{fig:stable}
\end{figure}
\end{Rem}

Note that in the Theorem the existence of the Feller process is assumed. It is clearly desirable to find conditions for the family of symbols $q(x,\xi)$ which already ensure the existence of the limit. This is part of ongoing research, for a survey see for example Jacob and Schilling \cite{JacoSchi2001}. It is remarkable that the conditions given in the usual constructions are much stronger then the requirements of the theorem above. Thus all these processes can be approximated by the Euler scheme. One obvious idea to find further sufficient conditions on the symbol would be to translate the conditions on the coefficients for example given by Stroock \cite{Stro2003} into conditions on $q(x,\xi)$ but so far no general criterion for this is known, compare \cite{Stro2003} 3.2.2. and Tsuchiya \cite{Tsuc1992}. In this context note that \eqref{assumptionboundedness} reflects the assumption of bounded coefficients in the SDE setting.

In the next section we give the necessary definitions and in Section \ref{proof} the proof of the Theorem is presented.

\section{Preliminaries} \label{def}
Let $C_\infty(\bbr^d)$ and $C_c^\infty(\bbr^d)$ be the continuous functions vanishing at infinity and the arbitrary often differentiable functions with compact support respectively. $\overline{B_0(1)}$ is the closed unitball in $\bbr^d$ and we use the notation $1_{\{\abs{g(y)}\leq 1 \}}$ for $1_{\overline{B_0(1)}}(g(y))$.

We consider $\R^d$ valued Markov processes in the
sense of Blumenthal and Getoor \cite{blumenthalget} and denote such
a process by $\textbf{X}=(\Omega, \cf, (\cf_t)_{t\geq 0},
(X_t)_{t\geq 0}, \bbp^x)_{x\in\bbr^d}.$ The expectation with respect to $\bbp^x$ is
denoted by $\bbe^x$.

A stochastic process $\textbf{X}$ is called \emph{Feller process} if the family of operators $(T_t)_{t\geq 0}$ defined by
\[
    T_t u(x) := \bbe^x(u(X_t)), \quad u\in C_\infty(\bbr^d),
\]
is a \emph{Feller semigroup}.  This is a strongly continuous contraction semigroup on $C_\infty(\bbr^d)$ which is positivity preserving. 
The semigroup and thus the corresponding process is called
conservative, if $T_t1 = 1$.

The infinitesimal generator $(A,\Dskript(A))$ of a Feller semigroup is defined by
\[
    A u := \lim_{t\to 0} \frac{T_t u - u}{t}\ \ \text{ on }\ \
    \Dskript(A) := \left\{u\in C_\infty(\bbr^d) \::\: \lim_{t\to 0} \frac{T_t u - u}{t} \text{\ exists strongly}\right\}.
\]
If the test functions $C_c^\infty(\bbr^d)$ are contained in the domain of the generator $A$ of a Feller semigroup, Courr\`ege \cite{Cour1964} showed that $A\big|_{C_c^\infty(\bbr^d)}$ has a representation as pseudo differential operator:
\begin{equation*}
    Au(x) = -q(x,D)u(x) = -\int_{\rn} q(x,\xi) e^{ix'\xi}\hat u(\xi)\,d\xi, \quad u\in C_c^\infty(\bbr^d)
\end{equation*}
where $\hat u(\xi) := (2\pi)^{-d}\int_\rn e^{-ix'\xi}\,u(x)\,dx$ is the Fourier transform of $u$. The function $q:\rn\times\rn\to\C$ is called the symbol of the pseudo differential operator. It is measurable and locally bounded in $(x,\xi)$. Furthermore it is continuous and negative definite (in the sense of Schoenberg) as a function of $\xi$, that is $\xi\mapsto q(x,\xi)$ admits for each $x\in\bbr^d$ the following L\'evy-Khinchine representation:
\begin{equation}\label{lkf}
q(x,\xi)=c(x)-i\ell(x)'\xi+\frac{1}{2} \xi'Q(x) \xi - \int_{y\neq 0} \Big( e^{i\xi'y}-1-i\xi'y 1_{\overline{B_1(0)}}(y) \Big) \ N(x,dy)
\end{equation}
where $c(x)\geq 0$, $\ell(x)\in\rn$, $Q(x)\in \R^{\dimension\times \dimension}$ is
positive semidefinite and $N(x,dy)$ is a kernel on
$\rn\times\mathcal{B}\left(\R^\dimension\setminus\{0\}\right)$ such
that $\int_{\rn\setminus\{0\}} (\norm{y}^2\wedge 1)\ N(x,dy)
<\infty$ for all $x$.
Especially one has
\begin{equation} \label{symbolrealpart}
\Re q(x,\xi) = c(x) + \frac{1}{2} \xi'Q(x) \xi + \int_{y\neq 0} \left(1 -\cos(\xi'y)\right)\ N(x,dy) \geq 0.
\end{equation}
In the following we will assume without loss of generality that every Feller process we encounter has c\`adl\`ag paths (cf.\,\cite{revu1999} Theorem III.2.7). The space of all c\`adl\`ag functions from $[0,\infty)$ to $\R^d$ is denoted by $D([0,\infty),\R^d)$ and convergence in this space is meant with respect to the Skorokhod $J_1$-topology (cf.\,\cite{Kall1997}).

A stochastic process is a L\'evy process if it has stationary and
independent increments and c\`adl\`ag paths. In particular note that
every L\,evy process is a Feller process with a symbol not depending
on $x$, i.e. $q(x,\xi)=\psi(\xi)$ and $\psi$ has a L\'evy-Khinchine
representation, see Sato \cite{Sato99} for further details.

Further we set $\widetilde{\cp}:=\cp\otimes \cb(\bbr^d)$ where $\cp$ is the predictable $\sigma$-algebra and $\cb(\bbr^d)$ denotes the Borel sets of $\bbr^d$. An integral with respect to a vector of processes is meant as matrix-vector multiplication, i.e. for a $d\times d$-matrix valued process $Y$ and and a $d$-dimensional vector valued process $X$ we write
\[
\int_0^t Y_{s-} \ dX_s = \left( \begin{array}{c} \sum_{j=1}^d \int_0^t Y^{1j} \ dX^{(j)}_s \\ \vdots \\ \sum_{j=1}^d \int_0^t Y^{dj} \ dX^{(j)}_s \end{array} \right).
\]
Integrals with respect to random measures are denoted by $H*\mu$ (cf.\,\cite{JacoShir2002} Section 2.1) and this integral is meant componentwise if $H$ is a vector.

If random variables $X$ and $Y$ are equal in distribution we write $X\stackrel{D}{=}Y$. Finally for a measure $N$ and a measurable function $g$ the image measure (push forward) is denoted by $N(g(\pushforward)\in dy)$.

\section{Proof of the Theorem} \label{proof}

Let $(\Omega, \cf, (\cf_t)_{t\geq 0}, (X_t)_{t\geq 0}, \bbp^x)_{x\in\bbr^d}$ be a c\`adl\`ag Feller process with generator $A$ such that $C_c^\infty(\bbr^d) \subset D(A)$.
The process admits the symbol $q:\bbr^d\times \bbr^d \to \bbc$ given by \eqref{lkf} with $c(x) = 0$ for all $x$, since \eqref{assumption-nokill} holds. 

Note that $A(C_c^\infty(\R^d)) \subset C_\infty(\R^d)$ and thus \eqref{assumption-nokill} implies by Schilling \cite{Schi1998} Theorem 4.4. that
\begin{equation} \label{assumptionfinely}
x \mapsto q(x,\xi) \text{\ \ is continuous for all $\xi$.}
\end{equation}
Furthermore $X_t$ is conservative by Theorem 5.2. in \cite{Schi1998} using \eqref{assumptioncore}-\eqref{assumption-nokill}.

Now the proof will be divided into three parts. First a result about Markov chain approximation of Feller processes is recalled and afterwards the SDE of a Feller process is calculated explicitly. Finally the characteristic functions of the increments of the Euler scheme to the SDE are calculated and the Markov chain approximation result is applied.

Given the assumptions the main theorem of B\"ottcher and Schilling \cite{BoetSchi2009} is applicable. This shows that
$$
    Y^{\stepsize} ([\cdot\,{\textstyle\frac{1}{\stepsize}}]) \xrightarrow{\stepsize \to 0} X_. \ \ \text{ in }\ D([0,\infty),\R^d)
$$
where $\left(Y^{\stepsize}(k)\right)_{k\in\N}$ is for each $\stepsize>0$ a Markov chain with initial value $Y^{\stepsize}(0):= X_0$ and transition kernel $\mu_{x,\stepsize}(dy)$ defined by
\begin{equation*}\label{charfun}
    \int_\rn e^{iy'\xi} \mu_{x,\stepsize}(dy) = e^{ix'\xi-\stepsize q(x,\xi)}.
\end{equation*}
To make our argumentation more self contained we note that this can
also be seen in the following way: let $U_{\stepsize}$ be the
transition operator corresponding to the kernel
$\mu_{x,\stepsize}(dy)$. By \eqref{assumptioncore} we can apply
Theorem 17.28 of Kallenberg \cite{Kall1997}. Using the mean value
theorem twice with suitable intermediate values $r,s \in
(0,\stepsize)$ and applying \eqref{assumptionboundedness} and
\eqref{symbolrealpart} it follows that
$$
   \left| \frac{e^{-\stepsize q(x,\xi)} - 1}{\stepsize} +  q(x,\xi)\right| = \left|-q(x,\xi) \left(e^{-sq(x,\xi)} - 1\right)\right| = \left|s q(x,\xi)^2 e^{-r q(x,\xi)}\right| \leq c^2 \stepsize \,(1+\left|\xi\right|^2)^2.
$$
Thus for $f\in C_c^\infty(\bbr^d)$ we obtain
$$    \left| \frac{U_{\stepsize} f(x) - f(x)}{\stepsize} - A f(x) \right|
    = \left| \int e^{ix'\xi} \left( \frac{e^{-\stepsize q(x,\xi)} - 1}{\stepsize} + q(x,\xi) \right) \hat f(\xi)\,d\xi\right| \leq c^2 \stepsize \int (1+|\xi|^2)^2 |\hat f(\xi)| \,d\xi \rightarrow 0
$$
for $\stepsize \downarrow 0$. This convergence is uniformly in $x$ and thus the $i)\Rightarrow iv)$ part of Theorem 17.25 \cite{Kall1997} implies that the Markov chain approximates the Feller process in $D([0,\infty),\R^d).$

In the following we will show that the approximation defined in \eqref{euler} for the SDE corresponding to $X_t$ coincides with the Markov chain in distribution, i.e.
\begin{equation}\label{eulercharfun}
\bar{X}_{(m+1)\cdot \stepsize} \text{ has the characteristic function }e^{ix'\xi} e^{-\stepsize q(x,\xi)} \text{ given that } \bar{X}_{m\cdot \stepsize} = x.
\end{equation}

For this we have first to find the SDE corresponding to $X_t$ explicitly.

By Schnurr \cite{Schn2009} Theorem 3.14 (see also Schilling \cite{Schi1998c}), using the assumption of conservativeness, one obtains that $(X_t)_{t\geq 0}$ is an It\^o process in the sense of Cinlar et.\,al.\,\cite{vierleute}, i.e. it is a strong Markov process which is a semimartingale with respect to every $\bbp^x$ and its characteristics are
\begin{align*}
B_t^{(j)}(\omega) &= \int_0^t \ell^{(j)} (X_s(\omega)) \ ds & j\in \{1,...,d\}, \\
C_t^{jk}(\omega)  &= \int_0^t Q^{jk}(X_s(\omega)) \ ds & j,k\in \{1,...,d\},    \\
\nu(\omega;ds,dy) &= N(X_s(\omega),dy) \ ds,
\end{align*}
with respect to the truncation function $h(y):=y1_{\overline{B_1(0)}}(y)$. By Cinlar and Jacod \cite{cinlarjacod81} Theorem 3.33 we obtain that on a suitable enlargement of the stochastic basis, the so called Markov extension, the process $(X_t)_{t\geq 0}$ is the solution of the following SDE. Let the Markov extension be
\[
(\widetilde{\Omega}, \widetilde{\cf}, (\widetilde{\cf}_t)_{t\geq 0}, (X_t)_{t\geq 0}, \widetilde{\bbp}^x)_{x\in\bbr^d}.
\]
On this space we have
\begin{equation}
\label{sde-feller}
\begin{split}
X_t=x &+\int_0^t \ell(X_{s-}) \ ds + \int_0^t \sigma(X_{s-}) \ d\widetilde{W}_s \\
   &+\int_0^t\int_{z\neq 0} k(X_{s-},z) 1_{\{\abs{k(X_{s-},z)}\leq1\}} \ \Big(\widetilde{\mu}(\cdot;ds,dz) - ds \widetilde{N}(dz)\Big) \\
   &+\int_0^t\int_{z\neq 0} k(X_{s-},z) 1_{\{\abs{k(X_{s-},z)}>1\}} \ \widetilde{\mu}(\cdot;ds,dz)
\end{split}
\end{equation}
where $\widetilde{W}$ is a $d$-dimensional Brownian motion, $\widetilde{\mu}$ is a Poisson random measure on $\roi{0}{\infty}\times \bbr\backslash \{0\}$ with dual predictable projection $dt \widetilde{N}(dz)$. Furthermore $\ell:\bbr^d \to \bbr^d, \ \sigma: \bbr^d \to \bbr^{d\times d}$ and $k:\bbr^d\times \bbr\backslash \{0\} \to \bbr^d$ are measurable functions and such that
\[
\widetilde{N}(k(X_s(\omega),\pushforward)\in
dy)ds=\nu(\omega;ds,dy)
\]
$\widetilde{\bbp}^x$-a.s. for every $x\in\bbr^d$ on the Markov extension (compare in this context Cinlar and Jacod \cite{cinlarjacod81} (3.9) and their remark following Theorem 3.7). Note that
\begin{align*}
\int_0^t\int_{z\neq 0} k(X_{s-},z)& \left(1_{\{\abs{k(X_{s-},z)}\leq1\}}-1_{\{|z|\leq 1\}}\right)\ \widetilde{N}(dz) ds\\
&=\bigg[\int_0^t\int_{0<|z|\leq 1} k(X_{s-},z) \ \Big(\widetilde{\mu}(\cdot;ds,dz) - ds \widetilde{N}(dz)\Big) +\int_0^t\int_{|z|> 1} k(X_{s-},z) \ \widetilde{\mu}(\cdot;ds,dz)\bigg]\\
&\ \ \ -\bigg[\int_0^t\int_{z\neq 0} k(X_{s-},z) 1_{\{\abs{k(X_{s-},z)}\leq1\}} \ \Big(\widetilde{\mu}(\cdot;ds,dz) - ds \widetilde{N}(dz)\Big)\\
&\ \ \ \ \ \ \  +\int_0^t\int_{z\neq 0} k(X_{s-},z) 1_{\{\abs{k(X_{s-},z)}>1\}}) \
    \widetilde{\mu}(\cdot;ds,dz)\bigg].\\
\end{align*}
The integral on the left hand side exists since representation \eqref{sde-feller} is valid and therefore either $k(X_{s-},z) \xrightarrow{z \to 0} 0$ or $\tilde{N}$ integrates constants at the origin.
Thus by a change of the cutoff function \eqref{sde-feller} is the same as
\begin{align*}
X_t=x &+\int_0^t \tilde{\ell}(X_{s-}) \ ds + \int_0^t \sigma(X_{s-}) \ d\widetilde{W}_s \\
   &+\int_0^t\int_{0<|z|\leq 1} k(X_{s-},z) \ \Big(\widetilde{\mu}(\cdot;ds,dz) - ds \widetilde{N}(dz)\Big) \\
   &+\int_0^t\int_{|z|> 1} k(X_{s-},z) \ \widetilde{\mu}(\cdot;ds,dz)
\end{align*}
where $\tilde{\ell}(x) = \ell(x) - \int_{z\neq 0} k(x,z) \left(1_{\{\abs{k(x,z)}\leq1\}}-1_{\{|z|\leq 1\}}\right)\ \widetilde{N}(dz)$. Thus $(X_t)_{t\geq 0}$ is the solution of an SDE of form \eqref{sde}.

In \cite{Schn2009} Theorem 5.7 (see also \cite{SchiSchn09}) it is shown that, given \eqref{assumptionfinely}, the It\^o process $(X_t)_{t\geq 0}$ has the symbol $p:\bbr^d\times\bbr^d\to \bbc$ given by
\begin{align} \label{symbolofeq}
p(x,\xi)= -i\ell(x)'\xi + \frac{1}{2} \xi'\sigma(x) \sigma(x)' \xi - \int_{z\neq 0} \Big(e^{i k(x,z)'\xi} -1 - ik(x,z)'\xi\cdot 1_{\overline{B_1(0)}}(k(x,z))\Big) \ \widetilde{N}(dz).
\end{align}
i.e.
\[
p(x,\xi):=- \lim_{t\downarrow 0}\bbe^x \frac{e^{i(X_{t \land T}-x)'\xi}-1}{t}
\]
for every first exit time $T$ of a compact set containing $x$.
The symbol $p(x,\xi)$ and the symbol $q(x,\xi)$ coincide for Feller processes by Corollary 4.5 of \cite{Schn2009}. For every fixed $x\in\bbr^d$ both functions are continuous and negative definite in the co-variable $\xi$. Since the L\'evy triplet of such a function is unique for a fixed cut-off function (cf.\,\cite{BergFors75} Theorem 10.8), we obtain $Q(x)=\sigma(x)\sigma(x)'$ and
\begin{eqnarray} \label{compensators}
N(x,dy)=\widetilde{N}(k(x,\pushforward)\in dy).
\end{eqnarray}
Now we define for fixed $t\geq 0$ and $x\in\bbr^d$ a process $(Y_{\stepsize})_{\stepsize \geq 0}$ by
\begin{align} \label{processy}
Y_{\stepsize} =x &+\int_t^{t+ \stepsize} \ell(x) \ ds + \int_t^{t+ \stepsize} \sigma(x) \ d\widetilde{W}_s \nonumber \\
   &+\int_t^{t+ \stepsize}\int_{z\neq 0} k(x,z) 1_{\{\abs{k(x,z)}\leq1\}} \ \Big(\widetilde{\mu}(\cdot;ds,dz) - ds \widetilde{N}(dz)\Big) \\
   &+\int_t^{t+ \stepsize}\int_{z\neq 0} k(x,z) 1_{\{\abs{k(x,z)}>1\}} \ \widetilde{\mu}(\cdot;ds,dz) \nonumber.
\end{align}

\begin{Lem}
For the process $Y:=(Y_{\stepsize})_{\stepsize\geq 0}$ defined above we obtain:
\begin{itemize}
  \item[a)] $Y$ is a L\'evy process.
  \item[b)] The following identity holds in distribution
                \begin{align*}
                  Y_{\stepsize}\stackrel{D}{=}x&+\stepsize \ell(x) + \sigma(x) W_{\stepsize} \\
                   &+\int_0^{\stepsize} \int_{y\neq 0} y \cdot 1_{\{\abs{y} \leq 1 \}} \ \Big( \mutilde(\cdot;ds,k(x,\pushforward)\in dy) - ds N(x,dy) \Big) \\
                   &+\int_0^{\stepsize} \int_{y\neq 0} y \cdot 1_{\{\abs{y}> 1 \}} \ \mutilde(\cdot;ds,k(x,\pushforward)\in dy)
                \end{align*}
  \item[c)] The characteristic function of $Y_{\stepsize}$ is
             \[
                \bbe^x\Big(e^{iY_{\stepsize}'\xi} \Big) = e^{ix'\xi} e^{-\stepsize q(x,\xi)}.
             \]
\end{itemize}
\end{Lem}

\textbf{Proof:} Fix $x\in\bbr^d$ and $t\geq 0$. The four integral terms in \eqref{processy} are stochastically independent and we will treat them separately. For the first two integrals all statements of the Lemma are easily obtained, since for every $t\geq 0$
\[
  \int_t^{t+ \stepsize} \ell(x) \ ds + \int_t^{t+ \stepsize} \sigma(x) \ d\widetilde{W}_s = \ell(x) \stepsize + \sigma(x) (\widetilde{W}_{t+\stepsize}-\widetilde{W}_t)
\]
where $\widetilde{W}_{t+\stepsize}-\widetilde{W}_t$ is again a Brownian motion having the same distribution as $(\widetilde{W}_{\stepsize})_{\stepsize \geq 0}$. For the integrals with respect to the (compensated) random measures we will have to proceed step-by-step.

a) First we show that the random measure $ds N(x,dy)$ is the dual predictable projection of the measure $\mutilde(\omega; ds, k(x,\pushforward)\in dy)$. To this end let $H:(\Omega,\bbr_+,\bbr^d)\to\bbr$ be positive and $\widetilde{\cp}$-measurable. Then we have
\begin{align*}
\bbe^x \Big(H(\cdot,s,y)*\mutilde(\cdot;ds,k(x,\pushforward)\in dy)\Big) &= \bbe^x\Big( H(\cdot,s,k(x,z))* \ \widetilde{\mu}(\cdot;ds,dz) \Big) \\
  &= \bbe^x\Big( H(\cdot,s,k(x,z)) * \ ds \widetilde{N}(dz) \Big) \\
  &= \bbe^x\Big( H(\cdot,s,y)* \ ds \widetilde{N}(k(x,\pushforward)\in dy) \Big) \\
  &= \bbe^x\Big( H(\cdot,s,y)* \ ds N(x,dy) \Big)
\end{align*}
where we have used that $ds \widetilde{N}(dz)$ is the dual predictable projection of $\widetilde{\mu}(\cdot;ds,dz)$ and \eqref{compensators} for the last equality.
By Theorem II.1.8 of \cite{JacoShir2002} we obtain that $ds N(x,dy)$ is the dual predictable projection of $\mutilde(\omega; ds, k(x,\pushforward)\in dy)$.

Now we are in the position to deal with the third integral term:
\begin{align*}
\int_t^{t+\stepsize} & \int_{z \neq 0} k(x,z) 1_{\{\abs{k(x,z)}\leq 1 \}} \ \Big( \mutilde(\cdot;ds,dz)-ds \widetilde{N}(dz) \Big) \\
&= \lim_{n\to\infty} \int_t^{t+\stepsize} \int_{\loi{-\infty}{-1/n}\cup\roi{1/n}{\infty}} k(x,z) 1_{\{\abs{k(x,z)}\leq 1 \}} \ \Big( \mutilde(\cdot;ds,dz)-ds \widetilde{N}(dz) \Big) \\
&= \lim_{n\to\infty} \int_t^{t+\stepsize} \int_{k(x,\cdot)(\loi{-\infty}{-1/n}\cup\roi{1/n}{\infty})} y 1_{\{\abs{y}\leq 1 \}} \ \Big( \mutilde(\cdot;ds,k(x,\pushforward)\in dy)-ds N(x,dy) \Big) \\
&= \int_t^{t+\stepsize} \int_{y\neq 0} y 1_{\{\abs{y}\leq 1 \}} \
\Big( \mutilde(\cdot;ds,k(x,\pushforward)\in dy)-ds N(x,dy) \Big).
\end{align*}
Let us emphasize that the integrals in the second and third line can be written as the difference of two integrals with respect to the respective random measures. Therefore these measures can be transformed one-by-one. The limit in third line exists and therefore does the limit in the second line, too.

In particular $\mutilde(\cdot;ds,k(x,\pushforward)\in dy)$ is a Poisson random measure by the structure of its compensator, because $N(x,dy)$ is a L\'evy measure for every fixed $x\in\bbr^d$.

The fourth term can now be written as
\[
\int_t^{t+\stepsize}\int_{z\neq 0} k(x,z) 1_{\{\abs{k(x,z)}>1\}} \
\mutilde(\cdot;ds,dz) = \int_t^{t+\stepsize}\int_{y\neq 0} y
1_{\{\abs{y}>1\}} \ \mutilde(\cdot;ds,k(x,\pushforward)\in dy).
\]
Putting the four terms together we obtain a L\'evy-It\^o decomposition (cf.\,Chapter 4 of \cite{Sato99}) of the process $Y=(Y_{\stepsize})_{\stepsize\geq 0}$, although the third and fourth integral are still `shifted'. In particular $Y$ is a L\'evy process.

b) It is enough to give the proof for the case $d=1$, because the integrals with respect to the (compensated) random measures are defined componentwise. This time we start with the fourth term of \eqref{processy}. Since $\mutilde(\cdot;ds,k(x,\pushforward)\in dy)$ is a Poisson random measure we have
\[
  \mutilde(\cdot;\loi{a}{b},k(x,\pushforward)\in C) \ \stackrel{D}{=} \  \mutilde(\cdot;\loi{t+a}{t+b},k(x,\pushforward)\in C)
\]
for $a<b$ and $C\in \cb(\bbr^1\backslash \{0\})$. Therefore we obtain
\[
  \int_0^\infty \int_{\{y>1\}} 1_{\loi{a}{b}}(s) \cdot 1_C(y) \ \mutilde(\cdot;ds,k(x,\pushforward)\in dy) \stackrel{D}{=}
  \int_0^\infty \int_{\{y>1\}} 1_{\loi{t+a}{t+b}}(s) \cdot 1_C(y) \ \mutilde(\cdot;ds,k(x,\pushforward)\in dy).
\]
In order to keep notation simple, we set
\[
  I(a,b;1_C):= \int_0^\infty \int_{\{y>1\}} 1_{\loi{a}{b}}(s) \cdot 1_C(y) \ \mutilde(\cdot;ds,k(x,\pushforward)\in dy).
\]
Now let $\varphi:\bbr^d\backslash \{0\}\to \bbr_+$ be a simple function, i.e. $\varphi$ can be written as $\varphi(y)=\sum_{j=1}^m d_j\cdot 1_{C_j}(y)$ where $m\in\bbn$ the $C_j$ are disjoint sets in $\bbr^d \backslash\{0\}$ and $d_j\geq 0$ for every $j\in \{1,...,m\}$. Since the $C_j$ are disjoint the random variables $(I(a,b,1_{C_i}))_{i=1,\ldots,m}$ and $(I(t+a,t+b,1_{C_i}))_{i=1,\ldots,m}$ are independent respectively. And thus we obtain
\[
  I(a,b;\varphi) \stackrel{D}{=} I(t+a,t+b;\varphi).
\]
Furthermore there exists a sequence of simple functions $(\varphi_n)_{n\in\bbn}$ such that $\varphi_n \uparrow \id$ on $(1,\infty)$. In the limit we obtain
\[
  I(a,b;\id) \stackrel{D}{=} I(t+a,t+b;\id)
\]
and analogous on $(-\infty,-1)$. The fact
\[
\Big(\loi{a}{b}\times \{y>1\} \Big) \cap \Big(\loi{a}{b}\times \{y<-1\} \Big)=\emptyset
\]
implies that the random variables
\[
\int_a^b \int_{\{y >1\}} y \ \mutilde(\cdot;ds,k(x,\pushforward)\in dy) \text{ and } \int_a^b \int_{\{y <-1\}} y \ \mutilde(\cdot;ds,k(x,\pushforward)\in dy)
\]
are independent. Therefore we obtain the representation of the
fourth term of b).

To deal with the third term we write
\begin{align*}
 \int_t^{t+\stepsize}& \int_{y\neq 0} y 1_{\{\abs{y}\leq 1 \}} \ \Big( \mutilde(\cdot;ds,k(x,\pushforward)\in dy)-ds N(x,dy) \Big) \\
&= \lim_{n\to\infty} \int_t^{t+\stepsize} \int_{\{\abs{y}>1/n\}} y 1_{\{\abs{y}\leq 1 \}} \ \Big( \mutilde(\cdot;ds,k(x,\pushforward)\in dy)-ds N(x,dy) \Big)
\end{align*}
and note that by using the same arguments as for the fourth term we
get
\[
\int_t^{t+\stepsize} \int_{\{\abs{y}>1/n\}} y 1_{\{\abs{y}\leq 1 \}} \ \mutilde(\cdot;ds,k(x,\pushforward)\in dy)=\int_0^{\stepsize} \int_{\{\abs{y}>1/n\}} y 1_{\{\abs{y}\leq 1 \}} \ \mutilde(\cdot;ds,k(x,\pushforward)\in dy).
\]
The integrals with respect to $ds N(x,dy)$ are clearly invariant with respect to a shift of $t$ because of the product structure and the non-randomness. Thus we obtain
\begin{align} \label{test}
Y_{\stepsize}\stackrel{D}{=}x&+\stepsize \ell(x) + \sigma(x) \widetilde{W}_{\stepsize} \nonumber \\
   &+\int_0^{\stepsize} \int_{y\neq 0} y \cdot 1_{\{\abs{y} \leq 1 \}} \ \Big( \mutilde(\cdot;ds,k(x,\pushforward)\in dy) - ds N(x,dy) \Big) \\
   &+\int_0^{\stepsize} \int_{y\neq 0} y \cdot 1_{\{\abs{y}> 1 \}} \ \mutilde(\cdot;ds,k(x,\pushforward)\in dy). \nonumber
\end{align}

c) The right-hand side of \eqref{test} is the L\'evy-It\^o decomposition of a L\'evy process. This process has the same one-dimensional distributions as $Y$. Therefore the characteristic functions of the two processes coincide and we obtain
\[
\bbe^x \Big(e^{iY_{\stepsize}'\xi} \Big) =e^{ix\xi} e^{-\stepsize q(x,\xi)},
\]
which proves the Lemma. \hfill $\square$

For $t = m\stepsize$ the Lemma shows that \eqref{eulercharfun} holds and thus the Theorem is proven.

\begin{ack}
We would like to thank an anonymous referee for carefully reading
the manuscript and offering useful suggestions which helped to
improve the paper.
\end{ack}



\end{document}